\documentclass[leqno,12pt]{article}

\usepackage[utf8]{inputenc}
\usepackage[T1]{fontenc}
\usepackage[english]{babel}
\usepackage[intlimits]{amsmath}
\usepackage{amsfonts}
\usepackage{amssymb}
\usepackage{bbm}

\usepackage{textcomp}
\usepackage{array}
\usepackage{enumerate}

\usepackage[amsmath, thmmarks]{ntheorem}
\usepackage{multicol}
\usepackage{units}
\usepackage{nicefrac}

\newcounter{theoremcounter}

\numberwithin{equation}{section}
\newtheorem{globaltheoremcounter}{theoremcounter}
\newtheorem{theoremcounter}{theoremcounter}[section]

\theoremstyle{plain}
\theoremseparator{.}

\theoremstyle{plain}
\newtheorem{globaltheorem}[globaltheoremcounter]{Theorem}

\newtheorem{lemma}[theoremcounter]{Lemma}
\newtheorem{theorem}[theoremcounter]{Theorem}

\newtheorem{proposition}[theoremcounter]{Proposition}
\newtheorem{remark}{Remark}

\theorembodyfont{\normalfont}
\theoremsymbol{\ensuremath{\Box}}
\newtheorem{example}{Example}
\newtheorem*{eproof}{Proof}

\newcommand{\N}{\mathbb{N}}
\newcommand{\Z}{\mathbb{Z}}
\newcommand{\Q}{\mathbb{Q}}

\newcommand{\diag}{\mathrm{diag}}
\newcommand{\tr}{\mathrm{tr}}

\newcommand{\Mat}[1]{\mathrm{M}_{#1}}
\newcommand{\GL}[1]{\mathrm{GL}_{#1}}

\newcommand{\RGL}[1]{\widetilde{\mathrm{GL}}_{#1}}
\newcommand{\Sp}[1]{\mathrm{Sp}_{#1}}
\newcommand{\RSp}[1]{\widetilde{\mathrm{Sp}}_{#1}}
\newcommand{\Inv}[1]{\mathrm{Inv}_{#1}}
\newcommand{\RInv}[1]{\mathrm{Inv}_{#1} ^{\mathrm{Res}}}

\newcommand{\compl}[2]{#1_{\mathfrak{#2}}}

\newcommand{\adelic}[1]{#1_{\mathbb{A}}}

\newcommand{\K}{\mathbb{K}}
\newcommand{\p}{\mathfrak{p}}
\newcommand{\oK}{\mathfrak{o}_\K}

\begin{document}
\title{Elementary divisor theory for the modular group over quadratic field extensions and quaternion algebras}
\author{Martin Raum\\\small{Martin.Raum@matha.rwth-aachen.de}\vspace{1.5em}\\ Lehrstuhl A f\"ur Mathematik\thanks{the author's present address is MPI für Mathematik, Vivatsgasse 7, 53111 Bonn, Germany}\\RWTH Aachen University\\52056 Aachen, Germany}
\date{\today}
\maketitle

\begin{abstract}
We develop an elementary divisor theory for the unimodular and the modular group over quadratic field extensions and quaternion algebras.
In particular, we investigate which sets of elementary divisors can occur.
Under an additional hypothesis we establish a correspondence of unimodular and modular double cosets.
\end{abstract}

\section{Introduction and statement of the results}

Elementary divisor theory has proven to be a useful tool to investigate many different structures. Namely, we consider two matrix groups $G_1,\,G_2$. For matrices $M,\,M'$ in a suitable monoid, we say they are equivalent, $M \sim M'$, if and only if there are $U \in G_1,\,V \in G_2$ such that $M' = U M V$. We want to analyze equivalence classes with respect to this relation.
For the simplest base rings $\Lambda$, which are principal ideal domains, and $G_1 = G_2 = \GL{n} (\Lambda)$ elementary divisor theory provides a classification for finitely generated modules.

This paper is mainly dedicated to $G_1 = G_2 = \Sp{n} (\Lambda)$ over maximal orders $\Lambda$ of either a quadratic field extension or a quaternion algebra $\Omega$ over some number field $\K$. We will call this the modular case. The related equivalence relation is $M \sim_m M' \Leftrightarrow \exists U,\,V \in \Sp{n}(\Lambda) \,:\,U M V = M'$. The group $\Sp{n}(\Lambda)$ will be defined below.

We will refer to $G_1 = G_2 = \GL{n} (\Lambda)$ as the unimodular case.
In this case the equivalence relation is $M \sim_u M' \Leftrightarrow \exists U,\,V \in \GL{n}(\Lambda) \,:\,U M V = M'$.

The idea to analyze $\sim$ is to find $U,\,V$ such that $U M V$ is a diagonal matrix.
Its entries are called elementary divisors. In case this is not possible, one tries to define suitable substitutes.
For example $\left(\begin{smallmatrix}3 & \rho \\ \rho & 3\end{smallmatrix}\right)$ with $\rho = \sqrt{-6}$ admits no diagonal form with respect to $G_1 = G_2 = \GL{2}(\Z [\rho])$. Nevertheless, it is useful to say its elementary divisors are the ideals $(3,\,\rho)$ and $(3)$. See Theorem \ref{trivialunimodulargenus} and \ref{trivialmodulargenus}.

With this general notion of elementary divisors in mind another concrete question is the following. Let $\Lambda \mid \Z$ be a quaternion maximal order ramified at $17$ and $\infty$. Is there a matrix with elementary divisors $\mathfrak{m}_2$ and $2$, where $\mathfrak{m}_2$ is a maximal left ideal over $2$? The answer is actually no, as we will see in Theorem \ref{unimodularexistence} and Theorem \ref{existenceofmodularmatrices}.
\vspace{1.5em}

To state the main results let $\K$ be a number field. Let $\Omega \mid \K$ be either a quadratic field extension or a quaternion algebra over $\K$ which is an integral domain. 
Fix a maximal order in $\Omega$ and denote it by $\Lambda$. The maximal order in $\K$ will be denoted by $\oK$.
Let $\adelic{\Omega}$ be the adeles of $\Omega$. Let $\adelic{\Lambda} \subseteq \adelic{\Omega}$ be the adeles which are integral at all finite places of $\Lambda$.
By convention we have $\Omega \subseteq \adelic{\Omega}$.
We will denote the completions at a finite place $\p$ by $\Omega_\p \subseteq \adelic{\Omega}$ and $\Lambda_\p$, respectively.

Fix $n \in \N$. In the quaternion case we have to assume $n \ge 2$.

The integral invertible matrices $M \in \Inv{n} (\Lambda) = \GL{n}(\Omega) \cap \Mat{n}(\Lambda)$ will be of interest in the unimodular case. In section \ref{sec:unimodularcase} we will define local elementary divisors $e_{\p,\,i} (M) \in \Lambda_\p$ for $i \in \{1,\ldots,n\}$ and finite places $\p$ of $\Lambda$.

We turn to the modular case.
The similitude matrices form the basic object of interest.
There is a non-trivial involution $\imath$ of $\Omega$ over $\K$.
In the quaternion case it restricts to the Galois conjugation on maximal subfields of $\Omega$.
We define $M^* = \imath(M)^{\mathrm{T}}$ where $\cdot\,^{\mathrm{T}}$ is the matrix transposition.
Set $J = \left(\begin{smallmatrix}& I_n \\ -I_n &\end{smallmatrix}\right)$ where $I_n \in \Mat{n}(\Lambda)$ is the identity matrix.
Given $m \in \mathrm{Fix}_{\imath} (\adelic{\Omega}) \setminus \{0\}$ we set
\begin{align*}
  \Delta_n (\adelic{\Omega},\, m)
=
  \{M \in \Mat{2 n}(\adelic{\Omega}) \,:\, M^* J M = m J\}
.
\end{align*}
This is the monoid of all similitude matrices.
If $R \subseteq \adelic{\Omega}$ is a subring and $m \in \mathrm{Fix}_{\imath}(R) \setminus \{0\}$ we will denote the $R$-valued similitude matrices by $\Delta_n (R,\,m) = \Delta_n (\adelic{\Omega},\, m) \cap \Mat{2 n}(R)$.
In particular, we will consider $\Delta_n (\adelic{\Lambda},\,m) = \Delta_n (\adelic{\Omega},\, m) \cap \Mat{2 n}(\adelic{\Lambda})$ and $\Delta_n (\Lambda,\,m) = \Delta_n (\adelic{\Omega},\, m) \cap \Mat{2 n}(\Lambda)$.

The extended symplectic group is
$
  \Sp{n}(\adelic{\Lambda})
=
  \bigcup_{m \in H} \Delta_{n}(\adelic{\Lambda},\, m)
$. Here $H = \mathrm{Fix}_{\imath}(\adelic{\Lambda}^\times)$ is the group of units of $\adelic{\Lambda}$ stable under $\imath$.
We set $\Sp{n}(\Lambda) = \Sp{n}(\adelic{\Lambda}) \cap \Mat{2 n}(\Lambda)$.

The norm of $a \in \adelic{\Omega}$ will be denoted by $N(a) = a \imath(a)$.
The normalized exponential valuation of an adele or an ideal with respect to a finite place $\p$ of $\Lambda$ will be denoted by $\nu_\p$.
We denote the class group of finitely generated, locally free $\Lambda$-left modules by $\mathrm{Cl}(\Lambda)$.
Notice that in our setting its elements correspond to actual isomorphism classes, not only to stable isomorphism classes.

\begin{globaltheorem}
\label{modularunimodularcorrespondence}
We assume that $\mathrm{Cl} (\oK)$ is trivial. Let $m \in \oK \setminus \{0\}$.
\nopagebreak

The equivalence classes $\Delta_n (\Lambda,\,m) / \sim_m$ and
\begin{align*}
\{M \in \Inv{n}(\Lambda) \,:\, \nu_\p (N(e_n(M))) \le \nu_\p (m) \text{ for all } \p \} / \sim_u
\end{align*}
are in one-to-one correspondence.
Here $\p$ runs through all finite places of $\Lambda$ and $e_n (M) := \prod_\p e_{\p,\,n} (M) \in \adelic{\Lambda}$.
\end{globaltheorem}

\begin{remark}
One interpretation of the theorem above is that the loss of degrees of freedom due to the symplectic structure of a module is compensated for by the dimension $2 n$.
It roughly behaves like a module without additional structures of dimension $n$.
\end{remark}

A more precise description of these equivalence classes will be provided in \ref{trivialunimodulargenus}, \ref{unimodularexistence} and \ref{existenceofmodularmatrices}.
Using the second main result we come to a more detailed understanding.

\begin{globaltheorem}
\label{trivialmodulargenus}
Let $M,\,M' \in \Delta_n (\Lambda,\, m)$ for some $m \in \oK \setminus \{0\}$.
There are $U,\, V \in \Sp{n}(\Lambda)$ such that $U M V = M'$ if and only if there are $\hat U,\, \hat V \in \Sp{n}(\adelic{\Lambda})$ such that $\hat U M \hat V = M'$.
\end{globaltheorem}

\begin{remark}
\label{modularimplementation}
The proof of this theorem as well as of Theorem \ref{trivialunimodulargenus} actually provides a way to construct $U,\,V$ from $\hat U,\,\hat V$ by the chinese remainder theorem. 
This might be used to efficiently enumerate right cosets in a double coset $\Sp{n}(\Lambda) M \Sp{n}(\Lambda)$.
Naturally, one can also decide whether $M \sim_m M'$ by a semimodular algorithm, which promises to be quite fast (cf. \cite{Lue02}).
Currently, there is only a rough implementation of the unimodular analog over number fields in Sage (\cite{sage}).

Notice that for small $n$ in the quaternion case it might be more efficients to use a splitting extension of $\Omega$ (cf. Example \ref{ex:quaternionsovernumberfield}).
\end{remark}

For every finite place $\p$ of $\Lambda$ we will derive local modular elementary divisors $e_{\p,\,i} \in \Lambda_\p$ for $i \in \{1,\ldots,n\}$ and finite places $\p$ of $\Lambda$.
We will see (cf. Proposition \ref{quaternionicmodularelementarydivisors}) that the local modular elementary divisors completely determine the \mbox{$\sim_m$-equivalence} class for integral similitude matrices.
Naturally, the question arises which sets of modular elementary divisors may occur.
To answer this question notice that every $e_{\p,\,i} \in \Lambda_\p$ corresponds to a $\Lambda$-left ideal.
Hence, it also defines an element of $\mathrm{Cl}(\Lambda)$.

\begin{globaltheorem}
\label{existenceofmodularmatrices}
We assume that $\mathrm{Cl} (\oK)$ is trivial.
\nopagebreak

We consider $m \in \oK \setminus \{0\}$ and $M \in \Delta_n (\adelic{\Lambda},\,m)$.
There is a similitude matrix $M' \in \Delta_n (\Lambda,\,m)$ with $e_{\p,\,i} (M) = e_{\p,\,i} (M')$ for all finite places $\p$ of $\Lambda$ and all $i \in \{1,\ldots,n\}$  if and only if $\sum_{\p,\,i} e_{\p,\,i} (M) = 0 \in \mathrm{Cl}(\Lambda)$.
\end{globaltheorem}
\begin{remark}
As we will see in Proposition \ref{quaternionicmodularelementarydivisors} the \mbox{$n$-th} modular elementary divisor of $M \in  \Delta_n (\Lambda,\, m)$ satisfies $\nu_\p ( N(e_n (M))) \le \nu_\p (m)$ for all finite places $\p$ of $\Lambda$.
\end{remark}

\vspace{1.5em}
We comment on related works and our hypotheses.

For general base rings and groups the question on elementary divisors turned out to be tough.
A particular example is Nakayamas's question.
In \cite{Na38} he asked for a set of invariants of modules over a ring $\Lambda$.
This question led to Guralnick's investigations of the unimodular case in \cite{Gu88,Gu88_pid}.
We give a more concrete proof for Theorem \ref{trivialunimodulargenus} which corresponds to \cite[Lemma 3.8]{Gu88}.
Theorem \ref{unimodularexistence} extends his work by describing the elementary divisors which can occur for integral invertible matrices.
There has been an investigation of this question in the case $\K = \Q$ and $\Omega$ hermitian over $\K$. In this setting Theorem \ref{unimodularexistence} is a classical result by Franz (cf. \cite{Fr34}).

Shimura investigated lattices in \cite{Sh64}. In our language his setting includes the case $G_1 = \GL{n},\, G_2 = \mathrm{U}_n (\phi)$, where $\phi$ is a hermitian form.
Recently, the unimodular and the modular case have been considered in \cite{Eb08, Eb09} by means of elementary methods.

Typically, $\Lambda$ will be a hermitian field extension or a definite quaternion algebra and $G_1 = G_2$ will be related to some additional structure imposed on $\Lambda$-modules.
Interesting example comprise quadratic (in some cases also called hermitian) modules, which play a central roll for modular forms (for example confer \cite{Bo98}).
Another group of interesting examples is given by symplectic modules, which are related to Siegel modular forms via Hecke theory.
In \cite{Ra09} the author comes back to this particular application.

If $\Omega$ is non-commutative we have to assume $n \ge 2$ to satisfy the Drozd condition (cf. \cite[section 2.9]{Gu88}).
One could also use the Eichler condition, but for the sake of simplicity we will not.
Confer \cite{Vi76} for a treatment of the cancellation properties in definite quaternion algebras.

In Theorem \ref{modularunimodularcorrespondence} and \ref{existenceofmodularmatrices} we have to assume that $\oK$ is a principal ideal domain.
This assumption is needed to guarantee that there is a decomposition of hermitian matrices into coprime nominator and denominator (cf. Lemma \ref{coprimematrixfraction}).

Our notion of the simplectic group differs from the one Krieg used in \cite{Kr87_ha}.
He set $\Sp{n}(\Lambda) = \Delta_n (\Lambda,\,1)$.
In our setting this is not appropriate, since there might be infinitely many units in $\oK$ and it might be hard to find a multiplicatively closed set of representatives of $(\oK \setminus \{0\}) / \oK^\times$.
\vspace{1.5em}

The paper is organized as follows.
In section \ref{sec:unimodularcase} we will investigate the unimodular case. In section \ref{sec:modularcase} we will define modular elementary divisors and we will show that they uniquely determine the \mbox{$\sim_m$-equivalence} classes. The question which sets of modular elementary divisors can occcure will be treated in section \ref{sec:existenceofmodularmatrices}. In the last section we will present some examples illustrating the theory.

\subsubsection*{Acknowledgements}
\textit{The author is indebted to Aloys Krieg and Julia Hartmann for helpful comments on the first versions of this paper. He also thanks the referees for their comments.}


\section{Unimodular equivalence}
\label{sec:unimodularcase}

We first fix our setting. In this section we assume that $\Omega \mid \K$ is either a field extension or central simple algebra (cs algebra) which is an integral domain.
As in the introduction we fix a maximal order $\Lambda \subseteq \Omega$.
Notice that in general there does not exist an involution $\imath$ of $\Omega$ over $\K$.
But we do not need it in this section.

We explain how to associate elementary divisors to a matrix $M \in \Inv{n}(\Lambda) = \GL{n}(\Omega) \cap \Mat{n}(\Lambda)$.
In the number field case there is a global construction.
Namely, we consider the determinantal divisors $d_i = \sum_{I,J \in \mathrm{Pow}_i (\{1,\ldots,n\})} \Lambda \det (M_{I,J})$.
Here $\mathrm{Pow}_i$ denotes the set of all subsets of cardinality $i$.
The minor of $M$ with rows $I$ and columns $J$ is denoted by $M_{I,J}$.

The global elementary divisors are the integral ideals $e_i (M) = d_i / d_{i-1}$ for $i \in \{1,\ldots,n\}$.
Here we set $d_{-1} = \Lambda$.
For every finite place $\p$ of $\Lambda$ we can define local elementary divisors $e_{\p,\,i} (M) = \Lambda_\p e_i (M)$.
As adeles they coincide with the elementary divisors obtained by considering the Smith normal form of $M$ with respect to $\GL{n} (\Lambda_\p)$.

This is how we are going to define the local elementary divisors in the cs algebra case.
Namely, we consider the Smith normal form of $M$ with respect to the unimodular group over $\Lambda_\p$, which is possible by \cite{Gu88_pid}.
This yields a diagonal matrix $\diag(e_{\p,\,1}(M), \ldots, e_{\p,\,n}(M)) = U M V$ for $U,\,V \in \GL{n} (\Lambda_\p)$ with $e_{\p,\,i}(M) \| e_{\p,\,i+1}(M)$ for $i \in \{1,\ldots,n\}$.
Here a block diagonal matrix is denoted by $\diag(\cdots)$ and $a \| b \Leftrightarrow \Lambda_\p b \Lambda_\p \subseteq \Lambda_\p a \cap a \Lambda_\p$. We say that $a$ is a total divisor of $b$ in $\Lambda_\p$.
The elements $e_{\p,\,i}(M) \in \Lambda_\p$ are unique up to similarity.
Namely, they are well defined up to multiplication by units from the left and from the right.
We will call them the local elementary divisors of $M$.
Notice that every $e_{\p,\,i}(M)$ corresponds to a  $\Lambda$-left ideal and we will occasionally use this identification.
\vspace{1em}

The local maximal order $\Lambda_\p$ is a left and right principal ideal domain.
Moreover, $\GL{n}(\Lambda_\p)$ is generated by elementary row and column operations.
To see this, we revise the following facts from \cite{Re75}.
The completion $\Omega_\p$ is isomorphic to $\Mat{r}(\mathfrak{D})$ for some local division algebra $\mathfrak{D}$.
The group $\GL{n}(\Lambda_\p)$ is isomorphic to $\GL{n r} (\mathfrak{D})$.
In $\mathfrak{D}$ every left ideal is a two-sided ideal.
Every element of $\mathfrak{D}$ is integral if and only if its norm is integral.
Hence, the group $\GL{n r}(\mathfrak{D})$ is generated by elementary row and column operations and these correspond to elementary row and column operations in $\GL{n}(\Lambda_\p)$.

For technical reasons we will have to consider a subgroup $\RGL{n}(\Lambda_\p) \subseteq \GL{n}(\Lambda_\p)$ for $n \ge 2$.
It is generated by the matrices $I_n + a I_{i,j}$ with $i \ne j,\,a \in \Lambda_\p$.
Here the identity matrix is denoted by $I_{n}$ and $I_{i,j} \in \Mat{n}(\Lambda_\p)$ is the matrix with only zeros but a single one in the $i$-th row and $j$-th column.
Define $\RGL{n}(\adelic{\Lambda}) \subseteq \GL{n}(\adelic{\Lambda})$ to be the group of all matrices $M$ such that $M_\p \in \RGL{n}(\Lambda_\p)$ for all finite places $\p$ of $\Lambda$.

Set $\Inv{n}(\adelic{\Lambda}) = \GL{n}(\adelic{\Omega}) \cap \Mat{n}(\adelic{\Lambda})$.
For all $M \in \Inv{n}(\adelic{\Lambda})$ we can find $m \in \oK \setminus \{0\}$ such that $m M^{-1} \in \Mat{n}(\adelic{\Lambda})$.

We consider the Smith normal form of $M$ with respect $\GL{n}(\Lambda_\p)$.
Every matrix in $\GL{n}(\Lambda_\p)$ is a product of a diagonal matrix and a matrix in $\RGL{n}(\Lambda_\p)$.
Thus, $U M V$ has Smith normal form for some $U,\,V \in \RGL{n}(\Lambda_\p)$.
Notice that since $\Lambda_\p$ is local the Smith normal form can be made explicit.
This is important for an implementation.

To prove Theorem \ref{trivialunimodulargenus} we need the following lemma.
It tells us that we can approximate elements of $\RGL{n}(\adelic{\Lambda})$ up to arbitrary ideals in $\Lambda$.
\begin{lemma}
\label{unimodularapproximationlemma}
Let $n \ge 2$, fix a finite place $\p$ of $\Lambda$ and consider $U \in \RGL{n} (\compl{\Lambda}{p})$.
Then for any $m \in \oK \setminus \{0\}$ there is a matrix $U' \in \GL{n} (\Lambda)$ such that $U \equiv U' \mod m \Lambda_\p$ and $U' \equiv I_n \mod m \compl{\Lambda}{q}$ for any finite place $\mathfrak{q} \ne \mathfrak{p}$.
The same holds for $U \in \RGL{n} (\adelic{\Lambda})$.
\end{lemma}
\begin{eproof}
The adelic case follows from the local one.
It suffices to prove the lemma for a set of generators of $\RGL{n} (\compl{\Lambda}{p})$.
Since all generators of $\RGL{n}(\Lambda_\p)$ differ only by one row or column from the identity matrix we can restrict to $n = 2$.

Hence, we have to consider $U = \left(\begin{smallmatrix}1 & a \\ & 1\end{smallmatrix}\right)$ or  $U =  \left(\begin{smallmatrix}1 & \\ a & 1\end{smallmatrix}\right)$.
We have $\compl{\Lambda}{p} / m \compl{\Lambda}{p} \cong \Lambda / \mathfrak{p}^{\nu_\mathfrak{p} (m)} \Lambda$ and this module has a complement in $\Lambda / m \Lambda$.
Thus we can choose $a' \in \Lambda$ with $a' \equiv a \mod m\compl{\Lambda}{p}$ and $a' \equiv 0 \mod m\compl{\Lambda}{q}$ for all $\mathfrak{q} \mid m\Lambda,\, \mathfrak{q} \ne \mathfrak{p}$.
The matrix  $U' = \left(\begin{smallmatrix}1 & a' \\ & 1\end{smallmatrix}\right)$ or  $U' =  \left(\begin{smallmatrix}1 & \\ a' & 1\end{smallmatrix}\right)$ is the desired approximation.
\end{eproof}
\begin{remark}
\label{constructingapproximations}
Fixing a $\oK$ basis for $\Lambda$ the problem of finding $a'$ amounts to applying the chinese remainder theorem repeatedly. Hence, if $U$ is given in terms of generators, $U'$ can be constructed efficiently. This is important in view of Remark \ref{unimodularimplementation} and \ref{modularimplementation}.
\end{remark}

\begin{theorem}
\label{trivialunimodulargenus}
Let $M,\, M' \in \Inv{n}(\Lambda)$.
Then there are $U,\,V \in \GL{n}(\Lambda)$ such that $U M V = M'$ if and only if there are $\hat U,\, \hat V \in \GL{n}(\adelic{\Lambda})$ such that $\hat U M \hat V = M'$.
\end{theorem}

\begin{remark}
\label{unimodularimplementation}
Remark \ref{modularimplementation} applies also to the unimodular case.
It is actually easier to find the Smith normal form with respect $\GL{n} (\Lambda_\p)$.
\end{remark}

\begin{eproof}[Theorem \ref{trivialunimodulargenus}]
If $U$ and $V$ are given, simply set $\hat U = U$ and $\hat V = V$.

Conversely, suppose we know $\hat U,\, \hat V$.
We can restrict to $n \ge 2$. Otherwise, by the global assumptions $\Omega$ is a field and the theorem is well known.
We may assume that $m {M'} ^{-1} \in \Mat{n}(\Lambda)$ for some $m \in \oK$.
Set $\hat M = M,\,U = V = I_n$ and choose an ordering of the places $\mathfrak{p} \mid m$.
Then for every such place proceed as follows.
We have $\hat U_\p \hat M \hat V_\p = M' \mod m \Lambda_\p$.
We use Lemma \ref{unimodularapproximationlemma} to obtain approximations of $\hat U_\p,\, \hat V_\p$ with respect to $m$ and denote them by $U_\p,\,V_\p \in \GL{n}(\Lambda)$.
Then we replace $\hat M,\, U$ and $V$ by $U_\p \hat M V_\p,\, U_\p U$ and $V V_\p$, respectively.
Now, we have $\hat M \equiv M \mod m\Lambda_\p$.

Finally, we obtain $\hat M$ such that $\hat M {M'}^{-1},\,M' \hat M^{-1} \in \Mat{n} (\adelic{\Lambda})$.
Thus, $\hat M$ and $M'$ are equivalent with respect to $\GL{n}(\Lambda)$ and so are $M$ and $M'$.
\end{eproof}

We complete this section by investigating which elementary divisors can occur for $M \in \Inv{n}(\Lambda)$.
For $\Omega = \K$ this has already been derived in \cite{Fr34} using elementary methods.

Suppose that $\Omega$ is a number field.
Let $M$ and $N$ be \mbox{$\Lambda$-modules} such that $(e_\mathfrak{p,\,i} (M) \ne \Lambda_\p \text{ for some }i) \Rightarrow  (e_\mathfrak{\p,\,i} (N) = \Lambda_\p \text{ for all }i)$ and $(e_\mathfrak{p,\,i} (N) \ne \Lambda_\p \text{ for some }i) \Rightarrow  (e_\mathfrak{\p,\,i} (M) = \Lambda_\p \text{ for all }i)$.
For every $\Lambda$-module $P$ which has elementary divisors $e_{\p,\,i} (P) = e_{\p,\,i} (M) \cdot e_{\p,\,i} (N)$ we have $M + N = P \in \mathrm{Cl} (\Lambda)$.
This also holds for $\Lambda$-left modules in the cs algebra case.
To see this apply \cite[Corollary 3.7]{Gu88} to $M \oplus N$.
More precisely, we have $M + N = M \oplus N \in \mathrm{Cl}(\Lambda)$ by definition and $P = M \oplus N \in \mathrm{Cl}(\Lambda)$ since genera are trivial.

With this preparation we can prove the following theorem, which tells us which sets of elementary divisors can occur for integral invertible matrices.
\begin{theorem}
\label{unimodularexistence}
Suppose $M \in \Inv{n}(\adelic{\Lambda})$.
Then there is a matrix $M' \in \Mat{n}(\Lambda)$ satisfying $e_{\p,\,i}(M) = e_{\p,\,i}(M')$ for all finite places $\p$ of $\Lambda$ and all $i \in \{1,\ldots,n\}$ if and only if $\sum_{\p,\,i} e_{\p,\, i} (M) = 0 \in \mathrm{Cl}(\Lambda)$.
\end{theorem}

\begin{remark}
In the number field case the last condition reflects the fact that the determinant of a matrix generates a principal ideal.
\end{remark}

\begin{eproof}[Theorem \ref{unimodularexistence}]
If $M' \in \Inv{n}(\Lambda)$ as in the statement exists then $\sum_{\p,\,i} e_{\p,\, i} (M) = \sum_{\p,\,i} e_{\p,\, i} (M') = 0$, since $\Lambda^n M' \cong \Lambda^n$ as $\Lambda$-left modules.

Conversly, suppose $\sum_{\p,\,i} e_{\p,\,i} (M) = 0$.
The elementary divisors $e_{\p,\,i}(M)$ correspond to $\Lambda$-left ideals.
We consider generators $e_{\p,\,i,\,1},\ldots,e_{\p,\,i,\,k_{\p,i}}$ of these ideals $e_{\p,\,i} (M) \subseteq \Lambda$.
Next, we enumerate all places $\p_1,\ldots,\p_l$ with non trivial elementary divisors of $M$.
We set $k = k_{p_1, 1} + \cdots + k_{p_l,n}$.
The image of \mbox{$H := \Lambda^{k} \hat M \le \Omega^n$} in $\mathrm{Cl}(\Lambda)$ vanishes. Here
\begin{equation*}
  \hat{M}
=
  \left(\begin{matrix}
  e_{\p_1,1,\,1}   &  0        &         & 0  \\
  \vdots      &           &         &    \\
  e_{\p_1,1,\,k_{\p_1,1}} &  0        &         & 0  \\
  0           & e_{\p_1,2,\,1} & 0       & 0  \\
              &           & \ddots  &    \\
              &           &         &  e_{\p_l,n,\,k_{\p_l,n}}
%
    \end{matrix}\right)
.
\end{equation*}
Thus, we find a matrix $F \in \GL{n}(\Omega)$ such that $H F = \Lambda^n$. By \cite[Corollary 3.7]{Gu88} there are matrices $U \in \GL{k} (\Lambda)$ and $V \in \GL{n} (\Lambda)$ such that
\begin{align*}
  U \hat{M} F V
&=
  \left(\begin{matrix}
  I_n & 0^{(n, k - n)}
  \end{matrix}\right)^\tr
\quad
\text{and hence}
\\
  U \hat{M}
&=
  \left(\begin{matrix}
  M' & 0^{(n, k - n)}
  \end{matrix}\right)^\tr
\end{align*}
for a suitable matrix $M' \in \Mat{n} (\Lambda)$. This matrix $M'$ has the desired local elementary divisors and, thus, it is invertible.
\end{eproof}


\section{Elementary divisors with respect to the modular group}
\label{sec:modularcase}

In this section we start examining the modular case using the assumptions given in the introduction.
In particular, $\Omega \mid \K$ will either be a quadratic field extension or a quaternion algebra.
In the second case we assume $n \ge 2$.
We keep the notation introduced in the previous section.

We are going to define local modular elementary divisors which we will also denote by $e_{\p,\,i}$. There will be no confusion, whether we mean unimodular or modular elementary divisors.
We will see (cf. Proposition \ref{quaternionicmodularelementarydivisors}) that the local modular elementary divisors completely determine the $\sim_m$-equivalence class for integral similitude matrices.

We first introduce some new notation.
We define $\phi_i : \Mat{2} (\adelic{\Omega}) \rightarrow \Mat{2n} (\adelic{\Omega})$ to be the following map.
We assign to a matrix $\left(\begin{smallmatrix} a & b \\ c & d \end{smallmatrix}\right)$ the identity matrix with $a$ in its $(i,i)$th place, $b$ in its $(i,n+i)$th place, $c$ in its $(n+i,i)$th place and $d$ in its $(n+i,n+i)$th place. 
Moreover, we define $\psi : \GL{n} (\adelic{\Omega}) \rightarrow \Sp{n} (\adelic{\Omega})$, which maps $U$ to $\diag(U, (U^*)^{-1})$.

If $\p = \imath(\p)$ we set $\RSp{1} (\Lambda_\p) = \langle I_n + a I_{1,2},\, I_n + a I_{2,1} \,:\, a \in \mathrm{Fix}_{\imath} (\Lambda_\p)\rangle$. 
For $n \ge 2$ we denote the group generated by $\phi_i (\RSp{1} (\Lambda_\p)),\,i \in \{1,\ldots,n\}$ and $\psi(\RGL{n} (\Lambda_\p))$ by $\RSp{n} (\Lambda_\p) \subseteq \Sp{n} (\Lambda_\p)$.
The local modular group $\Sp{n}(\Lambda_\p)$ is generated by modular diagonal matrices and $\RSp{n} (\Lambda_\p)$.
This is well known in the number field case since $\Lambda_\p$ is local.
In the quaternion algebra case confer the proof of \cite[II.2.3]{Kr_halfspacequaternions}.
Notice that using Kriegs notation $J_{(2)} \in \RSp{1} (\Lambda_\p)$.

In analogy to \ref{unimodularapproximationlemma} we find
\begin{lemma}
\label{modularapproximationlemma}
Fix a finite place $\p$ of $\Lambda$ and suppose $\p = \imath(\p)$.
Let $U \in \RSp{n} (\compl{\Lambda}{p})$.
Then for any $m \in \oK \setminus \{0\}$ there is a matrix $U' \in \Sp{n} (\Lambda)$ such that $U \equiv U' \mod m \Lambda_\p$ and $U' \equiv I_n \mod m \compl{\Lambda}{q}$ for any finite place $\mathfrak{q} \ne \mathfrak{p}$.
\end{lemma}
\begin{eproof}
We can restrict to generators.
The group $\RSp{n}(\Lambda_\p)$ is generated by the embeddings $\phi_i$ and $\psi$.
These embeddings are compatible with the inclusion $\Lambda \subseteq \Lambda_\p$ as well as quotients of $\Lambda$ and $\Lambda_\p$.
Since $\mathrm{Fix}_{\imath}(\Lambda_\p) = (\oK)_\p$, the result follows from Lemma \ref{unimodularapproximationlemma}.
\end{eproof}

\begin{remark}
Suppose $U \in \Sp{n} (\adelic{\Lambda})$ is trivial at all finite places $\p \ne \imath(\p)$ of $\Lambda$. If $U_\p \in \RSp{n} (\Lambda_\p)$ for all places $\p = \imath(\p)$ of $\Lambda$ there is an according approximation.
\end{remark}

\subsection{The number field case}

Suppose $\Omega \mid K$ is a quadratic field extension.
We consider a finite place $\p$ of $\Lambda$.
The ring of local integers $\Lambda_\p$ is a principal ideal domain.

Firstly, suppose $\p$ is either ramified or inert.
Every matrix $M \in \Delta_n (\Lambda, m)$ yields a set of local unimodular elementary divisors $e_{\p,\,1} (M), \ldots, e_{\p,\,2 n} (M)$, which also determine the \mbox{$\sim_m$-equivalence} class of $M$ over $\Lambda_\p$.

Secondly, we consider the split places.
The Galois conjugation acts non trivially on split prime ideals in $\Lambda$. For the associated finite places $\mathfrak{p}$ denote by \mbox{$\RInv{2 n} (\compl{\Lambda}{p}, l) \subseteq \Inv{2 n} (\compl{\Lambda}{p}) \times \N_0$} the set of all pairs $(M, l)$ such that $M$ has level $\p^l$. Namely, all matrices $M$ such that the fractional ideal generated by the entries of $M^{-1}$ is contained in $\p^{-l}$.
The following lemma shows that we can naturally associate elementary divisors to similitude matrices over $\compl{\Lambda}{p} \oplus \Lambda_{\imath(\p)}$.

\begin{lemma}
Let $m \in \mathrm{Fix}_{\imath} (\Lambda_\p \oplus \Lambda_{\imath(\p)})$.
There is a one-to-one correspondence of $\Delta_n (\compl{\Lambda}{p} \oplus \Lambda_{\imath(\p)}, m)$ and $\RInv{2 n} (\compl{\Lambda}{p}, \nu_\p (m))$.
It is induced by $M \mapsto M_\p$.

The correspondence extends to a correspondence of double cosets with respect to $\Sp{n} (\compl{\Lambda}{p} \oplus \Lambda_{\imath(\p)})$ and $\GL{2 n} (\compl{\Lambda}{p})$, respectively. Here $\GL{2 n} (\compl{\Lambda}{p})$ acts trivially on the second component of $\RInv{2 n} (\compl{\Lambda}{p})$.
\end{lemma}
\begin{eproof}
We write $M = (M_\p,\,M_{\imath(\p)})$ for some $M \in \Delta_n (\compl{\Lambda}{p} \oplus \Lambda_{\imath(\p)}, m)$ with $m \in \mathrm{Fix}_{\imath} (\compl{\Lambda}{p} \oplus \Lambda_{\imath(\p)}) \setminus \{(0,0)\}$.
The statements follow if we show that $M$ is completely determined by is $\p$-adic part $M_\p$.

The equation $M^* J M = m J$ is equivalent to $M_{\imath(\p)} = m J^{-1} (M_\p ^{-1})^* J$.
Since $M_{\imath(\p)}$ has to be integral, $m M_\p ^{-1}$ also has to be integral.
This completes the proof.
\end{eproof}
\begin{remark}
\label{splitplaceelementarydivisors}
Since there is a Smith normal form with respect to $\GL{2 n}(\Lambda_\p)$ the correspondence gives a set of elementary divisors for $M \in \Delta_n (\compl{\Lambda}{p} \oplus \Lambda_{\imath(\p)}, m)$.

Let $e'_1 \mid \cdots \mid e'_{2 n}$ be the elementary divisors of $M_\p$ with respect to the $\GL{2 n} (\Lambda_\p)$.
Then $e_{\p,\,i}(M) = e'_i$ for $i \in \{1,\ldots,n\}$ and $e_{\imath(\p),\,i} = m \imath(e'_{2 n - i})^{-1}$ are the elementary divisors with respect to $\Sp{n} (\compl{\Lambda}{p} \oplus \Lambda_{\imath(\p)})$.

The condition $\nu_\p (N(e_{\p,\,n} (M) \cdot e_{\imath(\p),\,n} (M))) \le \nu_\p (m)$ translates to $\nu_\p (e'_n) \le \nu_\p (e'_{n+1} (M))$. The condition $\nu_\p (e'_{2 n}) \le \nu_\p (m)$ corresponds to the $\imath(\p)$-integrality of $M$.
\end{remark}

We have to extend the modular approximation lemma to the split case.
\begin{lemma}
\label{modularapproximationlemma2}
Fix a finite place $\p$ of $\Lambda$ and suppose $\p \ne \imath(\p)$.
Let $U \in \RGL{2 n} (\Lambda_\p)$.
Then for any $m \in \oK \setminus \{0\}$ there is a matrix $U' \in \Sp{n}(\Lambda)$ such that $U \equiv U' \mod m\Lambda_\p$ and $U' \equiv I_n \mod m \Lambda_{\mathfrak{q}}$ for any finite place $\mathfrak{q} \not \in \{\p,\,\imath(\p)\}$ of $\Lambda$.
\end{lemma}
\begin{eproof}
We consider the generators $I_n + a I_{i,j}$ of $\RGL{2 n}(\Lambda_\p)$ with $i \ne j$.
If $i,j \le n$ or $i,j > n$ we can find approximations along the lines of the proof of \ref{unimodularapproximationlemma}.
Suppose $i \le n,\, j > n$.
We have to consider two cases.
Firstly, suppose $j = i + n$.
We choose some $a' \in \Lambda$ such that $a' \equiv a \mod m\Lambda_\p$, $a' \equiv 1 \mod m\Lambda_{\imath(\p)}$ and $a' \equiv 0 \mod m\Lambda_{\mathfrak{q}}$ for all $\mathfrak{q} \not \in \{\p,\,\imath(\p)\}$.
The matrix $I_n + a' \imath(a') I_{i,j}$ is the desired approximation.
Secondly, suppose $j \ne i + n$.
We choose some $a' \in \Lambda$ such that $a' \equiv a \mod m\Lambda_\p$ and $a' \equiv 0 \mod m\Lambda_{\mathfrak{q}}$ for all $\mathfrak{q} \ne \p$.
Now, $I_n + a' I_{i,j} + \imath(a') I_{j-n,i+n}$ is a suitable matrix.

The case $i > n,\,j \le n$ is the same.
\end{eproof}

Finally, we define the adelic version of $\RSp{n}$.
Namely, $M \in \RSp{n}(\adelic{\Lambda})$ if and only if $M_\p \in \RSp{n}(\Lambda_\p)$ for all $\p = \imath(\p)$ and $M_\p \in \RGL{2 n}(\Lambda_\p)$ for all $\p \ne \imath(\p)$.
Matrices in $\Sp{n} (\adelic{\Lambda})$ are up to diagonal matrices elements of $\RSp{n}(\adelic{\Lambda})$.

\begin{eproof}[Theorem \ref{trivialmodulargenus}]
Proceed as in the proof of Theorem \ref{trivialunimodulargenus}.
Namely, given $M \in \Delta_n (\Lambda,\,m)$ and $ \hat U,\, \hat V \in \RSp{n}(\adelic{\Lambda})$ such that $M' = \hat U M \hat V$, repeatedly choose approximations of $U_\p$ and $V_\p$.
This will finally yield $\hat M = U M V$ with $U,\,V \in \Sp{n}(\Lambda)$ such that $\hat M (M')^{-1},\, M' \hat M^{-1} \in \Mat{2 n}(\adelic{\Lambda})$.
\end{eproof}

\subsection{The quaternionic case}
This case turns out to behave slightly different from the unimodular case. Namely, we will see that there are matrices which are equivalent with respect to the unimodular group but not with respect to the modular group.

In this subsection we assume that $\Omega$ is a quaternion algebra over $\K$. In particular, let $n \ge 2$.
The set of all finite places of $\Lambda$ is fixed elementwise by the automorphism $\imath$.

We first consider elementary divisors.
\begin{proposition}
\label{quaternionicmodularelementarydivisors}
Let $M \in \Delta_n (\Lambda_\p,\, m)$ for some $m \in \mathrm{Fix}_{\imath} (\Lambda_\p) \setminus \{0\}$.
Then there are matrices $U,\,V \in \Sp{n} (\compl{\Lambda}{p})$ such that $U M V = \diag(e_{\p,\,1},\ldots,e_{\p,\,n},m \imath(e_{\p,\,1})^{-1}, \ldots, m \imath(e_{\p,\,n})^{-1})$.

Here $\diag(e_{\p,\,1},\ldots,e_{\p,\,n})$ is in Smith normal form with respect to the unimodular group and $\nu_\mathfrak{p} (N(e_{\p,\,n}(M))) \le \nu_\p (m)$.
These modular elementary divisors are unique up to similarity.
Moreover $U$ and $V$ can always be chosen in $\RSp{n} (\compl{\Lambda}{p})$.

We  will call $e_{\p,\,1} (M),\ldots,e_{\p,\,n} (M)$ the local modular elementary divisors of $M$.
\end{proposition}
\begin{eproof}
The existence follows along the lines of \cite[II.2.2 and II.2.3]{Kr_halfspacequaternions}.

We will prove the uniqueness.
If $\Lambda_\p$ is ramified, it follows from the unimodular elementary divisor theorem (cf. \cite{Gu88_pid}). Since the modular and unimodular elementary divisors $e_{\p,\,i} (M)$ for $i \in \{1,\ldots,n\}$ coincide.

Hence, we assume that $\Lambda_\p$ splits.
We are going to prove the uniqueness in analogy to \cite[Theorem 6]{Kr87_ha}.

Let $M = \diag(e_1,\ldots,e_n,m  \imath(e_1)^{-1}, \ldots, m \imath(e_n)^{-1})$ and $M' = \diag(\diag(e'_1,\ldots,e'_n,m \imath(e'_1)^{-1}, \ldots, m \imath(e'_n)^{-1})$ be two modular Smith normal forms satisfying the restrictions given in the proposition.
Suppose $M = U M' V$ for $U,\,V \in \Sp{n}(\Lambda_\p)$.
We have to prove that $e_i$ and $e'_i$ are similar for all $i \in \{1,\ldots,n\}$.

The local quaternion algebra $\Lambda_\p$ is isomorphic to $\Mat{2}((\oK)_\p)$.
Let $\pi$ be a prime element of $(\oK)_\p$.
We can assume that $e_i = \diag(\pi^{\eta_{i,1}},\,\pi^{\eta_{i,2}})$ and $e_i = \diag(\pi^{\eta'_{i,1}},\,\pi^{\eta'_{i,2}})$ with $\eta_{i,1} \le \eta_{i,2}$ and $\eta'_{i,1} \le \eta'_{i,2}$.

We fix a useful notation.
For any $a \in \Lambda_\p$ we set $\epsilon(a) = \pi^l$ for the greatest $l \in \N$ such that $a / \pi^l \in \Lambda_\p$.
In particular, $\epsilon(e_i) = \pi^{\eta_{i,1}}$ and $\epsilon(e'_i) = \pi^{\eta'_{i,1}}$.

Due to the uniqueness of the Smith normal form with respect to $\GL{4 n} ((\oK)_\p)$ we know that $e_i = e'_i$ for $i \in \{1,\ldots,n-1\}$.
We decompose $V$ such that $U M \diag(W, s {W^*}^{-1}) = M' \hat V$ with $\hat V = \left(\begin{smallmatrix} A & B \\ C & D\end{smallmatrix}\right)$ and $A_{\bullet,n} = (0,\ldots,0,a)^{\mathrm{T}}$ for some $W \in \GL{n} (\Lambda_\p)$. Here $s \in \Lambda_{\p}^\times$ is the similitude of $V$ such that $\hat V$ has similitude $1$.

We will only consider the $n$-th column of the left and right hand side of the equation above.
Its right $\mathrm{gcd}$ will be denoted by $g$.
We have $\epsilon (g) \mid \epsilon(e_n) = \pi^{\eta_{n,1}}$.
To see this notice that if we write down the $n$-th column of $W$ as a $2 \times 2 n$ matrix at least two rows will contain elements of $(\oK)_\p ^\times$. Combine this with the fact that $\eta_{i,1},\eta_{i,2} \le \eta_{n,1}$ for all $i \in \{1,\ldots,n-1\}$.

Let $(c_1,\ldots,c_n)^{\mathrm{T}}$ be the $n$-th column of $C$.
We can deduce from the symplectic relation $M^* J M = J$ that $\imath(a) c_n = \imath(c_n) a$.
Hence, $a = \alpha v$ and $c_n = \gamma v$ for some $\alpha,\,\gamma \in \oK$ and some $v \in \Lambda_\p$.
We consider the $n$-th column of the right hand side.
It is equal to $(0,\ldots,0,e'_n \alpha v,m \imath(e_1)^{-1} c_1, \ldots,m \imath(e_{n-1})^{-1} c_{n-1}, m \imath(e'_n)^{-1} \gamma v)$.
We can choose a right divisor $d$ of $e'_n v$ which is similar to $e'_n$.
Since $e'_n$ is a total divisor of $m \imath(e_i)^{-1}$ for all $i \in \{1,\ldots,n-1\}$ as well as a right divisor of $m \imath{e'_n} ^{-1}$, we find that $d$ is a right divisor of $g$.

Summarizing these first steps, we have $\epsilon(e'_n) = \epsilon(d) \mid \epsilon(g) \mid \epsilon(e_n)$.
Due to symmetry we can deduce $\epsilon(e'_n) = \epsilon(e_n)$.

To finish the proof we consider the $n$-th column $(w_1,\ldots,w_n)^{\mathrm{T}}$ of $W$.
There are some $\hat x_i \in \Lambda_\p$ such that $\sum_i \hat x_i w_i = 1$.
Set $x_i = \hat x_i \diag(\pi^{\eta_{n,2} - \eta_{i,1}}, \pi^{\eta_{n,2} - \eta_{i,2}})$.
Then $\sum_i x_i e_i w_i = \pi^{\eta_{n,2}} = s' N(e_n) / \epsilon(e_n)$ for some $s' \in (\oK)_\p ^{\times}$.
We have chosen $g$ such that it is a right $\mathrm{gcd}$ of $U (e_1 w_1, \ldots, e_n w_n, 0,\ldots,0)^{\mathrm{T}}$.
Since $U$ is invertible it is a right $\mathrm{gcd}$ of $(e_1 w_1, \ldots, e_n w_n, 0,\ldots, 0)^{\mathrm{T}}$.
Hence $g$ is a divisors of $\pi^{\eta_{n,2}} \in \Lambda_\p$.
Since $d$ is a right divisor of $g$ and since $d$ is similar to $e'_n$, we can deduce that $e'_n$ is a divisor of $\pi^{\eta_{n,2}}$.
So, $\eta_{n,2} \ge \eta'_{n,2}$.
Due to symmetry the result follows.
\end{eproof}

\begin{remark}
If $\Lambda_\p$ is split it can occur that the modular elementary divisors do not satisfy $e_{\p,\,n} \mathop{\|} m \imath(e_{\p,\,1})^{-1}$.
Hence, some matrices are equivalent with respect to $\GL{2 n} (\compl{\Lambda}{p})$, but form distinct equivalence classes with respect to $\Sp{n} (\compl{\Lambda}{p})$.
This phenomenon was first described in \cite{Kr87_ha}.
\end{remark}

Since over $\Lambda_\p$ for fixed $m$ we established the uniqueness of the modular elementary divisors $e_{\p,\,1},\ldots,e_{\p,\,n}$ up to similarity, we can prove the most important result of this section follows along the lines of the preceding section.
\begin{eproof}[Theorem \ref{trivialmodulargenus}]
We proceed as in the number field case.
According to Proposition \ref{quaternionicmodularelementarydivisors} it suffices to apply Lemma \ref{modularapproximationlemma} repeatedly.
\end{eproof}


\section{Existence of modular matrices}
\label{sec:existenceofmodularmatrices}

In this section we have to assume that $\oK$ is a principal ideal domain. Otherwise we are not able to proof Lemma \ref{coprimematrixfraction}, that is essential.

We first define some generalizations of coprime hermitian pairs defined in \cite{Br49}.
We call a pair of matrices $(A, B)$ with $A,\,B \in \Mat{n} (\Lambda)$ satisfying $A B^* = B A^*$ a hermitian or quaternionic pair, respectively.
Such a pair is called coprime if and only if all local elementary divisors of the matrix $\left( \begin{smallmatrix}A & B\end{smallmatrix}\right)$ are trivial.
Furthermore two pairs $(A,\,B)$ and $(A',\,B')$ are associated to each other if and only if $A'B^* = B' A^*$.

A matrix satisfying $M = M^*$ will be called hermitian, also in the quaternionic case.

\subsection{The number fields case}\label{ssec:exhermitiancase}

Assume that $\Omega \mid \K$ is a quadratic field extension.
This section extends Braun's treatment of the hermitian modular group in \cite{Br49}. 
Notice that since we are assuming that $\oK$ is a principal ideal domain all proofs which are provided there and which we are going to use apply to any quadratic field extension $\Omega \mid \K$.

We are going to prove the first existence theorem.
\begin{eproof}[Theorem \ref{existenceofmodularmatrices}]
Suppose a matrix $M$ satisfies the given restrictions.
By Theorem \ref{unimodularexistence} we can simply choose an appropriate $N \in \Inv{n}(\Lambda)$ and set $M' = \diag(N, m (N^*)^{-1})$.

Conversely, it suffices to prove that any $M' \in \Delta_n (\Lambda,\,m)$ is \mbox{$\sim_m$-equivalent} to a block diagonal matrix $\diag(N, m (N^*)^{-1})$ such that $e_{\p,\,i} (N) = e_{\p,\,i} (M')$.

Firstly, set $M' = \left(\begin{smallmatrix}A & B \\ C & D\end{smallmatrix}\right)$. Considering the Smith normal form over $\adelic{\Lambda}$ and applying the approximation Lemmas \ref{modularapproximationlemma} and \ref{modularapproximationlemma2} we can assume that $M'$ has Smith normal form modulo $m^{n+1} \Lambda$.
In particular, $A \in \Inv{n}(\Lambda)$ since the determinant $A$ does not vanish modulo $m^{n+1}$.

Secondly, we can find a coprime hermitian pair $(E,\,-F)$ associated to $(A,\,C)$ in analogy to \cite[Theorem 1]{Br49}.
We may obtain a modular matrix $T = \left(\begin{smallmatrix}E & F \\ G & H\end{smallmatrix}\right)$ with suitable $G$ and $H$ in analogy to \cite[Lemma 1]{Br51}.
This yields
\begin{align*}
  T M'
=
  \left(\begin{matrix}
  E A + F C & E B + F D \\
            & G B + H D
  \end{matrix}\right)
.
\end{align*}
We set $\hat A = E A + F C,\, \hat B = E B + F D$ and $\hat D = G B + H D$.
Notice that $\hat A$ has the same elementary divisors as $A$.
We have
\begin{align*}
  T M'
=:
  \left(\begin{matrix}
  \hat A & \hat B \\
         & \hat D
  \end{matrix}\right)
=
  \left(\begin{matrix}
   I_n & \hat B {\hat D}^{-1} \\
       & I_n
  \end{matrix}\right)
  \left(\begin{matrix}
   \hat A & \\
    & \hat D
  \end{matrix}\right)
.
\end{align*}

We want to show that $\hat B {\hat D}^{-1} \in \Mat{n} (\Lambda)$.
Then the proof will be complete.
Since $T M' \in \Delta_n (\Lambda, m)$ we have 
\begin{align*}
  \hat B {\hat D}^{-1}
&
=
  m^{-1} \hat B {\hat A}^*
=
  m^{-1} (E B + F D) (A^* E^* + C^* F^*)
\\&
=
    E (m^{-1} B) (A^* E^* + C^* F^*)
  +  F (m^{-1} D A^*) E^*
  + F D (m^{-1} C^*) F^*
.
\end{align*}
This sum is integral since $M'$ has Smith normal form modulo $m \Lambda$.
\end{eproof}

\subsection{The quaternionic case}

In analogy with the hermitian case we want to show that for each quaternionic pair there is an associated quaternionic pair which is coprime.
In the next tree lemmas we proceed just as in Braun \cite{Br49} in the hermitian case.
\begin{lemma}\label{localhermitiannormalform}
Let $H = H^* \in \Mat{n} (\compl{\Lambda}{p})$.
If $\mathfrak{p} \nmid 2 \compl{\Lambda}{p}$ or $\compl{\Lambda}{p}$ splits, there is a unitary matrix $U \in \GL{n} (\compl{\Lambda}{p})$ such that $U H U^*$ has Smith normal form over $\compl{\Lambda}{p}$.
Otherwise, we can find a matrix $U \in \GL{n} (\compl{\Lambda}{p})$ such that $U H U^* = \diag(B_1, \ldots,B_r)$ for some $r \in \N$. 
For every $i \in \{1,\ldots,r\}$ we either have $B_i \in \compl{\Lambda}{p}$ or $B_i = \left(\begin{smallmatrix} h_{11} & \imath(h_{12}) \\ h_{12} & h_{22} \end{smallmatrix}\right) \in \Mat{2} (\compl{\Lambda}{p})$ satisfying $\nu_\mathfrak{p} (h_{11}),\, \nu_\mathfrak{p} (h_{22}) > \nu_\mathfrak{p} (h_{12})$.
\end{lemma}
\begin{eproof}
The ramified case can be proven like the hermitian one. If $\compl{\Lambda}{p}$ splits, we proceed as follows.
Since $H$ is hermitian the upper left entry of $H$ is in $(\oK)_\p$. Hence, by applying a suitable transformation $U$ and considering $\compl{\Lambda}{p} \cong \Mat{2}((\oK)_\p)$ we obtain
\begin{align*}
  U H U^*
=
  \left(\begin{matrix}
  \left(\begin{smallmatrix}
  \alpha p^m &            \\
             & \alpha p^m
  \end{smallmatrix}\right)
  & \cdots \\
  \left(\begin{smallmatrix}
   \beta_1 p^{n_1} &  \\
          & \beta_2 p^{n_2}
  \end{smallmatrix}\right)
  \\
  0 & \ddots \\
  \vdots & 
  \end{matrix}\right)
.
\end{align*}
Here $n_1 \le n_2$, $\alpha \in (\oK)_\p ^\times,\, \beta_1, \beta_2 \in \{0,1\}$ and $\beta_1 \ge \beta_2$. We can suppose that $m \le n_1$ since otherwise we apply
\begin{align*}
\left(\begin{matrix}
 \left(\begin{smallmatrix}
 1 & 0 \\
 0 & 1
 \end{smallmatrix}\right)
&
 \left(\begin{smallmatrix}
 1 & 0 \\
 0 & 0
 \end{smallmatrix}\right)
&
 0^{(1, n -2)}
\\
 0
&
 1
&
 0^{(1, n -2)}
\\
 0^{(n-2,1)}
&
 0^{(n-2,1)}
&
 I_{n-2}
\end{matrix}\right)
.
\end{align*}
Here, $0^{(i,j)}$ denotes the $i \times j$ zero matrix.
Now, we continue as in the ramified case to reduce the first column to $\left(\begin{smallmatrix}* & 0 \cdots 0\end{smallmatrix}\right)^{\mathrm{T}}$.

By first using induction on the number of non-reduced columns and then reordering the diagonal matrix we complete the proof.
\end{eproof}

\begin{lemma}
\label{coprimematrixfraction}
We assume that $\mathrm{Cl} (\oK)$ is trivial.
Suppose that $H \in \Inv{n}(\Omega)$ is hermitian. Then there are coprime matrices $N, M \in \Inv{n}(\Lambda)$ such that $H = N^{-1} M$.
\end{lemma}
\begin{eproof}
We choose $m \in \mathfrak{o}_K \setminus \{0\}$ such that $m H \in \Mat{n}(\Lambda)$ and consider its local elementary divisors $e_{\mathfrak{p},\,1},\ldots,e_{\mathfrak{p},\,n}$ for some finite place $\p$ of $\Lambda$.

First fix a finite place $\p$ of $\Lambda$.
Suppose $\p \nmid 2 \Lambda$ or $\Lambda_\p$ is split.
By the preceding lemma we have $e_{\mathfrak{p},\,i} = \compl{\Lambda}{p} e'_{\mathfrak{p},\,i} $ for some $e'_{i,\,\p} \in \oK$.
In the exceptional case, namely if $\p \mid 2 \Lambda$ and $\Lambda_\p$ is ramified, consider the blocks $B_i$ which are elements of $\Inv{2}(\Lambda_\p)$.
They have local elementary divisors $b_i,\, b_i$ for some $b_i \in \Lambda_\p$.
Since $\Lambda_\p$ is ramified $(\Lambda_\p b_i)^2$ corresponds to the principal ideal $\Lambda \p \subseteq \Lambda$.
Hence, there is a matrix in $\Inv{2}(\Lambda)$ with local elementary divisors $b_i,\,b_i$ over $\Lambda_\p$, which is integrally invertible over all other places of $\Lambda$ according to \ref{unimodularexistence}.

We choose $m' \in \oK$ such that $m' (mH)^{-1} \in \Mat{n}(\Lambda)$.
We consider an approximation such that $\hat U m H \hat U^* \equiv \diag(B_{\p,\,1},\ldots,B_{\p,\,r_\p}) \mod m' \Lambda_\p$ for every finite place $\p \mid m'$ of $\Lambda$.
As explained in the paragraph above we can find elements $\hat B_{\p,\,i} \in \oK$ or matrices $\hat B_{\p,\,i} \in \Inv{2} (\Lambda)$ which are have exactly the local elementary divisors of $B_{\p,\,i}$ over $\Lambda_\p$.
We can assume that $\hat B_{\p,\,i} \equiv 1 \mod m'\Lambda_{\mathfrak{q}}$ or $\hat B_{\p,\,i} \equiv \left(\begin{smallmatrix}1 & 0 \\ 0 & 1 \end{smallmatrix}\right) \mod m'\Lambda_{\mathfrak{q}}$ for all finite places $\mathfrak{q} \ne \p$ of $\Lambda$.

We build block diagonal matrices $\diag(\hat B_{\p,\,1}, \ldots, \hat B_{\p,\,r_\p})$ for each of these places $\p$.
Now, we can separate the blocks which correspond to elementary divisors $\nu_\p (e_{\mathfrak{p},\,i}) \le \nu_\p (m)$ of $m H$.
This yields block diagonal matrices $\hat N_\p$ such that $(\prod_\p \hat N_\p) \hat U H \hat U^* \in \Mat{n}(\Lambda)$. The order of the product is not important since every $N_\p$ is trivial modulo $m'\Lambda_{\mathfrak{q}}$ for $\p \ne \mathfrak{q}$.
Setting $N = (\prod_\p \hat N_\p)\hat U$ and $M = N H$ we are done.
\end{eproof}

Now, it is quite easy to derive the desired statement on quaternionic pairs.
\begin{lemma}\label{associatedhermitianpair}
For every quaternionic pair $(A,\, B)$ with $A \in \Inv{n} (\Lambda)$ there is an associated coprime quaternionic pair $(A',\,B')$.
\end{lemma}
\begin{eproof}
Since $(A,\, B)$ is a quaternionic pair, $H = A^{-1} B$ is hermitian. Choose a suitable matrix $L \in \Mat{n}(\Lambda)$ such that $H + L \in \Inv{n}(\Lambda)$ and choose coprime $A_0,\,B_0 \in \Inv{n}(\Lambda)$ such that $(H + L)^{-1} = A_0 ^{-1} B_0$.
The result follows by setting $A' = B_0$ and $B' = A_0 - B_0 L$.
\end{eproof}

Every coprime pair can be completed to a modular matrix. This can be proven in analogy with the hermitian case (cf. \cite[Lemma 1]{Br51}).
The main result of this section can now be deduced along the lines of the number field case.
\begin{eproof}[Theorem \ref{existenceofmodularmatrices}]
We assume that $M = \left(\begin{smallmatrix}A & B \\ C & D\end{smallmatrix}\right)$ has Smith normal form modulo $m^{n+1}$.
According to Lemma \ref{associatedhermitianpair} we can find a coprime quaternionic pair $(E,\,-F)$ associated to $(A,\,C)$.
The same calculations as in the proof of the number field case yield the result.
\end{eproof}

Finally we prove the relation between unimodular and modular equivalence classes. Here we consider the number field case as well as the quaternionic case.
\begin{eproof}[Theorem \ref{modularunimodularcorrespondence}]
By the proofs of Theorem \ref{existenceofmodularmatrices} every matrix $M \in \Delta_n (\Lambda,\,m)$ is equivalent to a block diagonal matrix $\diag(N,\,m (N^*)^{-1})$.
Here $e_{\p,\,i} (N) = e_{\p,\,i} (M)$ for all finite places $\p$ of $\Lambda$ and all $i \in \{1,\ldots,n\}$.
By Theorem \ref{trivialmodulargenus} and \ref{trivialunimodulargenus} they determine $M$ and $N$, respectively.
In the non-split number field case the restriction $\nu_\p(N(e_{n})) \le \nu_\p (m)$ on the elementary divisors is clearly satisfied.
For all split places it is satisfied as explained in Remark \ref{splitplaceelementarydivisors}.
In the quaternion case it is satisfied according to Proposition \ref{quaternionicmodularelementarydivisors}
\end{eproof}


\section{Examples}

We will consider some concrete examples in detail to illustrate the theory.

The first two examples illustrate the fact that there is no primary decomposition in the associated Hecke algebras (cf. \cite{Eb08, Ra09}).

\begin{example}[Decomposition of double cosets; hermitian case]
We set $\Omega = \Q(\rho),\, \rho = \sqrt{-6}$. The matrix $M = \left(\begin{smallmatrix} 2 & \rho \\ 4 & \rho \end{smallmatrix}\right)$ has elementary divisors $(2,\rho),\, (6, 2 \rho)$. Its double coset with respect to the modular group does not admit a decomposition into $2$- and $3$-components. This would only by possible if there were matrices with local elementary divisors $((2,\rho),(2))$ and $((1),(3, \rho))$. This is not possible as Theorem \ref{unimodularexistence} tells us.
According to Theorem \ref{modularunimodularcorrespondence} the same holds for the modular matrix $\diag(M, 12 (M^*)^{-1})$.
\end{example}

\begin{example}[Decomposition of double cosets; quaternionic case]
We can construct a similar quaternionic example.
The first try would be to find a quaternion algebra over $\Q$, which contains $\Q(\rho)$.
This would have to be ramified at $2$ and $3$. This quaternion algebra is indefinite. Hence, we have to consider another example.

Let $\Omega = \Q((-17,-3))$. Namely, we consider the quaternion algebra over $\Q$ with generators $i,j,k=ij = -j i$ satisfying $i^2 = -17,\,j^2 = -3$.
The maximal order $\Lambda$, that we will consider, is spanned by $1$, $h_1 = 1/2 + j/2$, $h_2 = 1/2 + i/2 + j/6 + k/6$ and $h_3 = -1/2 + j/6 - k/3$.
We have $N(1) = N(h_1) = 1$ and $N(h_2) = N(h_3) = 6$.
So $\mathfrak{m}_2 = {}_\Lambda(2, h_2)$ is a maximal ideal over $2$ and $\mathfrak{m}_3 = {}_\Lambda(3,h_2)$ is a maximal ideal over $3$.
Neither of them is a principal ideal.
We want to find a matrix with elementary divisors ($\mathfrak{m}_2,2)$ and $(1,\mathfrak{m}_3)$.
According to the proof of Theorem \ref{unimodularexistence} we have to consider
\begin{align*}
\left(\begin{matrix}
2   &   \\
h_2 &   \\
    & 6 \\
    & 2 h_2
\end{matrix}\right)
.
\end{align*}
Via elementary row and column operations from the left we find it is equivalent to
\begin{align*}
  M
&=
  \left(\begin{matrix}
  -2   & 6 \\
  -h_2 & 2 h_2
  \end{matrix}\right)
\end{align*}

Again, the modular analog is $\diag(M, 12 (M^*)^{-1})$.
\end{example}

\begin{example}[Quaternion algebra over a number field]
\label{ex:quaternionsovernumberfield}
We consider $\K = \Q(\rho)$ with $\rho = \sqrt{-6}$ and $\Omega = \K((-3,\,-17))$ which is not a left principal ideal domain.
The maximal order we will use is the restriction of the full integral matrix order in $\K(\eta) \otimes \Omega$ to $\Omega$.
Here, $\eta = \sqrt{-3}$ and $\Omega \hookrightarrow \K(\eta) \otimes \Omega$ via $i \mapsto \left(\begin{smallmatrix}\eta & \\ & \eta\end{smallmatrix}\right)$ and $j \mapsto \left(\begin{smallmatrix}&1 \\ -17 & \end{smallmatrix}\right)$.
We consider the matrix
\begin{align*}
  T
=
  \left(\begin{matrix}
  6 \rho          & 3 \rho ( 1 + i) \\
  3 - i + 3 j + k & 6 \rho
  \end{matrix}\right)
.
\end{align*}
and want to determine its elementary divisors.
We first need to find some $m \in \oK$ such that $m T^{-1} \in \Mat{2}(\Lambda)$. The minimal choice is $m = (4 \rho - 1) \cdot 9 \cdot 4 \rho$. Hence, we have to consider the places $(4 \rho - 1) \oK,\, 3\oK$ and $\rho \oK$.
The resulting local elementary divisors are $(1, \mathfrak{m}_{4 \rho - 1}),\,(\mathfrak{m}_3,3 \mathfrak{m}_3)$ and $(2 \mathfrak{m}_\rho, 2 \rho$).
Here $\mathfrak{m}_\p$ denotes a maximal left ideal over a place $\p$.
\end{example}


\bibliography{bibliography}{}
\bibliographystyle{plain}

\end{document}